\documentclass[suppldata]{interact}

\usepackage{epstopdf}
\usepackage{subfigure}
\usepackage{color}
\usepackage{natbib}
	\definecolor{ao(english)}{rgb}{0.0, 0.5, 0.0}

		\definecolor{auburn}{rgb}{0.43, 0.21, 0.1}
\theoremstyle{plain}

\theoremstyle{definition}

\theoremstyle{remark}

\usepackage{amssymb}
\usepackage{amsthm}
\usepackage{amsmath}

\usepackage{amsfonts}
\usepackage{subfigure}
\usepackage{enumerate}

\newcommand{\bqan}{\begin{eqnarray}}
	\newcommand{\eqan}{\end{eqnarray}}
\newcommand{\vmu}{\boldsymbol{\mu}}
\newcommand{\vSigma}{\boldsymbol{\Sigma}}

\newcommand{\vX}{\boldsymbol{X}}
\newcommand{\whsigma}{\widehat{\sigma}}

\newcommand{\bqa}{\begin{eqnarray*}}
	\newcommand{\eqa}{\end{eqnarray*}}
\numberwithin{equation}{section}
\theoremstyle{plain}
\newtheorem{thm}{Theorem}[section]

\begin{document}




\title{Multiplication-Combination Tests for Incomplete Paired Data}

\author{
	\name{Lubna Amro$^1$\thanks{CONTACT: Lubna Amro. Email: lubna.amro@uni-ulm.de}, Frank Konietschke$^2$\thanks{\hspace{48pt}Frank Konietschke. Email: fxk141230@utdallas.edu} and Markus Pauly$^1$\thanks{\hspace{48pt}Markus Pauly. Email: markus.pauly@uni-ulm.de}}
	\affil{$1$ Institute of Statistics, University of Ulm, Helmholtzstrasse 20, 89081 Ulm, Germany.}
	\affil{$2$ Department of Mathematical Sciences, The University of Texas at Dallas, 75080 Richardson,TX, U.S.A.}
}

\maketitle

\begin{abstract}
We consider statistical procedures for hypothesis testing of real valued functionals of matched pairs with missing values. In order to improve the accuracy of existing methods, we propose a novel multiplication combination procedure. 
 {\color{black} Dividing the observed data into dependent (completely observed) pairs and independent (incompletely observed) components, it is based on combining separate results of adequate tests for the two sub datasets. Our methods can be applied for parametric as well as semi- and nonparametric models and 
make efficient use of all available data. In particular, the approaches are flexible and can be used to test different hypotheses in various models of interest. This is exemplified by a detailed study of mean- as well as rank-based apporaches.} Extensive simulations show that the proposed procedures are more accurate than existing competitors. {\color{black} A real data set illustrates the application of the methods.}    
\end{abstract}


\section{Introduction and General Idea}\label{int}

 Repeated measure designs are widely employed in scientific and medical research. The simplest repeated measure design occurs when the subjects are observed repeatedly under two different treatments or time points \textendash called the matched pairs design. The main difficulty that arises when dealing with matched pairs is the problem of missing values. Since this problem occurs frequently in practice, we intend to develop a general test procedure that can handle matched pairs with missing values for different parametric as well as nonparametric models.
We consider inference methods for testing null hypotheses formulated in terms of real valued functionals $\theta=\theta(F_1,F_2)$ in matched pairs. 
Here, $F_1$ and $F_2$ are the marginal distribution functions of the i.i.d random vectors $ \mathbf{X}_j= (X_{1j}, X_{2j})',\, j=1,\ldots,n$. Specific examples of interest are: 
\begin{itemize}
 \item $\theta = \mu_1 - \mu_2$, the mean-difference functional (for $\mu_i = E(X_{i1}),i=1,2$) with null hypothesis $H_0^\mu:\{\theta=0\}$ (two-sided) or $H_0:\{\theta\leq 0\}$, or\\[-1ex]
 \item $\theta = p = P(X_{11}< X_{22}) + P(X_{11}= X_{22})/2$ the {\color{black}Wilcoxon-}Mann-Whitney effect with null hypothesis $H_0^p:\{p=1/2\}$ (two-sided) or $H_0:\{p\leq 1/2\}$.
\end{itemize}
For the case of complete observations, this issue has been analyzed in detail, see e.g. \cite{janssen1999nonparametric}, \cite{munzel1999nonparametric}, \cite{munzel2002exact, konietschke2012studentized}, and {\color{black}\cite{rietveld2017paired}.}
Here, we discuss the more involved problem, where some of the components may be missing{\color{black}, i.e. not all subjects were observed under both treatment conditions or time points}.  Several test procedures are developed in the literature to tackle this problem.  {\color{black}However, many of them are only developed under specific assumptions using parametric or semi-parametric mixed models, see for example \cite{lindsey1999models, diggle2002analysis, verbeke2009linear} for intensive textbook treatments. Typical underlying assumptions are symmetry or even bivariate normality which are hard to verify in practice. Moreover, these procedures are usually non-robust to deviations and may result in inaccurate decisions caused by possibly inflated or conservative type-I error rates, see e.g. the simulation studies in \cite{xu2012accurate, konietschke2012ranking, samawi2014notes, amropauly2017} and \cite{fong2017rank}. To overcome these problems the typical recommendation is to use 
 adequately weighted studentized test statistics for the underlying paired and unpaired two sample problem, see e.g. \cite{samawi2014notes} or \cite{amropauly2017} for the mean functional and \cite{gao2007nonparametric}, \cite{konietschke2012ranking} and \cite{fong2017rank} for nonparametric situations.}
  In case of small or moderate sample sizes, however, estimating large quantiles of the test statistics'  distributions is challenging, leading to possibly inflated type-I error rates. To overcome this gap, we propose a novel multiplication-combination approach: {\color{black} We sort the data into completely and incompletely observed random vectors}
\bqan \label{model: missing}
\underbrace{\binom{X_{11}^{(c)}}{X_{21}^{(c)}},\dots, \binom{X^{(c)}_{1n_c}}{X^{(c)}_{2n_c}}}_{\vX^{(c)}}, \underbrace{\binom{X_{11}^{(i)}}{--},\dots,\binom{X_{1n_1}^{(i)}}{--},\binom{--}{X_{21}^{(i)}},\dots,\binom{--}{X_{2n_2}^{(i)}}}_{\vX^{(i)}}.
\eqan 
{\color{black} Here, we throughout assume a MCAR} (missing completely at random) mechanism for simplicity, where $\vX^{(c)}$ is independent of $\vX^{(i)}$. Let   $n=n_c+n_1+n_2$ denote the total number of subjects and let $N=2n_c+n_1+n_2$ denote the total number of observations, where  $n_c$ is the number of complete cases, and $n_i$ is the number of incomplete {\color{black} observations in component $j\neq i, i,j\in\{1,2\}$.}

Let    
$\varphi^{(c)} {\color{black} = \varphi^{(c)}(\vX^{(c)})}$ and $\varphi^{(i)} {\color{black} = \varphi^{(i)}(\vX^{(i)})}$ denote adequate tests for the null hypothesis of interest that are computed upon $\vX^{(c)}$ and $\vX^{(i)}$ separately. Then $\varphi^{(c)}$ and $\varphi^{(i)}$ are also independent. Therefore, by calculating each test at level $\alpha^{1/2}$ individually, a {\it multiplication-combination test (MCT)} $\varphi = \varphi^{(c)}\cdot \varphi^{(i)}$  is achieved by multiplying the separate test results at significance level $\alpha\in(0,1)$. The advantage is that estimating $(1-\alpha^{1/2})$-quantiles {\color{black} (of the underlying test statistics)} is usually more accurate than the estimation of common $(1-\alpha)$-quantiles. Moreover, the whole idea even works if we calculate the independent tests at level $\alpha^{\gamma}$ and $\alpha^{1-\gamma}$ ($\gamma>0$), respectively. Here, $\gamma$ may even depend on the sample sizes $n_1,n_2$ and $n_c$. We directly formulate this in the following theorem, where $H_0$ may be any null hypothesis of interest.
\begin{thm}\label{theo 1} {\color{black}For $\alpha_1,\alpha_2\in(0,1)$ set $\alpha =\alpha_1\cdot\alpha_2$.}\\
	(i)  If $\varphi^{(c)}$ and $\varphi^{(i)}$ are level $\alpha_1$ and $\alpha_2$ tests for $H_0$, respectively, then 
	$\varphi=\varphi^{(c)}\cdot\varphi^{(i)}$ is a level $\alpha$ test for $H_0$. Moreover, if both are exact (i.e. fulfill $E_{H_0}(\varphi^{(c)})=\alpha_1$ and $E_{H_0}(\varphi^{(i)})=\alpha_2$), $\varphi$ is exact as well (i.e.  $E_{H_0}(\varphi)=\alpha$).\\
	(ii) If $\varphi^{(c)}$ and $\varphi^{(i)}$ are of asymptotic level $\alpha_1$ and $\alpha_2$ respectively 
	(i.e. $E(\varphi^{(c)})\to\alpha_1$ and $E(\varphi^{(i)})\to\alpha_2$ under $H_0$ as $\min(n_c,n_1,n_2)\to \infty$)), then $\varphi$ is of asymptotic level $\alpha$ as well.     
\end{thm}

Typical choices for $\alpha_1$ and $\alpha_2$ are $\alpha_1=\alpha^{\frac{n_c}{n}}$ and $\alpha_2=\alpha^{\frac{n_1+n_2}{n}}$ 
or $\alpha_1=\alpha^{\frac{2n_c}{N}}$ and $\alpha_2=\alpha^{\frac{n_1+n_2}{N}}$ which reduce to 
$\alpha_1=\alpha_2=\sqrt{\alpha}$ in 'balanced' cases  with $n_c = n_1+n_2$ or $2n_c = n_1+n_2$, respectively. 

We note that the independence of the tests $\varphi^{(c)}$ and $\varphi^{(i)}$ is essential for the validity of the MCT. Thus, depending on the choice of tests, even specific MAR (missing at random) situations can be tackled. In addition, the approach is applicable for different parametric as well as nonparametric models since only proper tests  are needed for the corresponding paired and unpaired two-sample designs. 

The remainder of the paper is organized as follows: In Section \ref{semi} we provide a new test procedure for mean comparisons in incomplete paired data under a semi-parametric framework. Rank-based test procedures for 
nonparametric models including hypotheses formulated in terms of the Wilcoxon-Mann-Whitney effect are proposed in Section \ref{nonparametric}. Afterwards, in Section \ref{alt}, alternative methods for matched pairs with missing values are reviewed to select 'competitors' for the subsequent simulation study conducted in Section \ref{sim}. A real data example is analyzed in Section \ref{example} and the paper closes with a discussion in Section \ref{dis}.

\section{Semi-parametric procedure} \label{semi}
In this section, we employ our suggested approach to develop a new procedure for testing the null hypothesis $H_0^{\mu}:{\{\mu_1=\mu_2\}}$ against the one-sided alternative ${\{\mu_1>\mu_2\}}$ or the two-sided alternative ${\{\mu_1\neq\mu_2\}}$. We consider Model (\ref{model: missing}) and assume that the first components $X_{1g}^{(c)},X_{1k}^{(i)}$ are i.i.d. with mean 
$\mu_1$ and variance $\sigma_1^2\in(0,\infty)$ and the second components
$X_{2g}^{(c)},X_{2\ell}^{(i)}$ are i.i.d. with mean $\mu_2$ and variance $\sigma_2^2\in(0,\infty)$ for $g=1,\dots,n_c, k=1,\dots,n_1, l=1,\dots,n_2$. 
Also, the complete pairs $(X_{1g}^{(c)},X_{2g}^{(c)})'$ are assumed to be i.i.d. with mean vector $\vmu=(\mu_1,\mu_2)'$ and an unstructured covariance matrix $\vSigma>0$ that allows for {\color{black} heteroscedastic variances.} 
\\
{\color{black}
To develop new procedures for testing $H_0^\mu: \{\mu_1=\mu_2\}$ with regard to our novel idea we only have to take two simple steps: Choose two independent test statistics that are adequate for the paired and unpaired setting, respectively, together with proper critical values.} For the completely observed pairs, we therefore suggest to use the well known paired $t$-type statistic
\bqan \label{T1}
T_{t}= T_{t}(\vX^{(c)}) = \frac{n_c^{-1}\sum_{g=1}^{n_c} D_g}{\sqrt{\whsigma^2_c/n_c}} =  \frac{\overline{D}_\cdot}{\sqrt{\whsigma^2_c/n_c}},
\eqan
where $D_g=X_{1g}^{(c)}-X_{2g}^{(c)}$ denote the differences of the first and second component for $g=1,\ldots,n_c$ and 
$\whsigma^2_c=(n_c-1)^{-1}\sum_{g=1}^{n_c} (D_g-\overline{D}_\cdot)^2$ is their empirical variance. {\color{black} Under the null hypothesis and additionally assuming normality, the statistic $T_t$ follows a t-distribution with $n_c-1$ degrees of freedom ($T(n_c-1)$). 
Moreover, due to the CLT, the distribution of $T_t$ can be approximated by a ($T(n_c-1)$)-distribution in our general (possibly non-normal) set-up. Thus,} 
  asymptotic level $\alpha_1$-tests are given by $\varphi_{t}=1\{T_t>t_{n_{c}-1,1-\alpha_1}\}$ and $\varphi_{t,2}=1\{\left|T_t\right|>t_{n_{c}-1,1-\alpha_1/2}\}$ 
 for the one- and two-sided case, respectively. Here, $t_{f,u}$ denotes the $u$-quantile of the t-distribution with $f$ degrees of freedom.

 However, the $T(n_c-1)$-approximation may become poor in case of small samples or skewed distributions, leading to possibly liberal or conservative test decisions, see e.g. \cite{KP14}. To overcome this problem and to enhance its small sample properties, a permutation version of the paired $t-$test has been recommended in \citep{janssen1999testing, KP14}. It is based on randomly changing the component for each pair and was shown to be exact if the underlying distribution is invariant under permuting the two components (e.g. for $0$-symmetric differences) and asymptotically exact in general. The corresponding one-sided permutation test is denoted by 
 $\varphi_{t,p}=1\{T_t>c_{t,p}(1-\alpha_1)\}$, where $c_{t,p}(1-\alpha_1)$ is the $(1-\alpha_1)$-quantile of the permutation distribution of $T_t$.

{\color{black} Next,} for the incomplete vectors, we suggest to apply the Welch-type test statistic \nocite{welch1938significance,welch1947generalization}
\bqan \label{T2}
T_{w}= T_{w}(\vX^{(i)}) = \frac{\overline{X}_{1\cdot}^{(i)} - \overline{X}_{2\cdot}^{(i)}}{\sqrt{\whsigma_1^2/n_1 + \whsigma_2^2/n_2}},
\eqan
where $\overline{X}_{1\cdot}^{(i)} = n_1^{-1}\sum_{k=1}^{n_1} X_{1k}^{(i)}$ and $\overline{X}_{2\cdot}^{(i)} = n_2^{-1}\sum_{\ell=1}^{n_2} X_{2\ell}^{(i)}$ are the sample means and 
$\whsigma_1^2=  (n_1-1)^{-1}\sum_{k=1}^{n_1} (X_{1k}^{(i)} - \overline{X}_{1\cdot}^{(i)})^2$ and 
$\whsigma_2^2 =  (n_2-1)^{-1}\sum_{\ell=1}^{n_2} (X_{2\ell}^{(i)} - \overline{X}_{2\cdot}^{(i)})^2$ the corresponding empirical variances. An asymptotic level $\alpha_2$-test for the one-sided case is given by $\varphi_{w}=1\{T_{w}>t_{\hat{v},1-\alpha_2}\}$  
and for the two-sided case by $\varphi_{w,2}=1\{\left|T_{w}\right|>t_{\hat{v},1-\alpha_2/2}\}$, where the degrees of freedom $\hat{v}$ are calculated by Welch's method.

 These Welch-tests are, however, known to be rather conservative procedures if the underlying distribution is non-symmetric, see e.g. \cite{janssen1997, janssen2005monte}. Thus, to enhance their power \cite{janssen1997, janssen2005} proposed a studentized permutation test in the Welch-type statistic ($T_{welch}$) that is based on randomly permuting the pooled sample, see also \cite{chung2013, pauly2015asymptotic}. 
 Similar to $\varphi_{t,p}$ the resulting Welch permutation test $\varphi_{w,p} = 1\{T_t>c_{w,p}(1-\alpha_2)\}$ (in the one-sided case) is an exact procedures if the data is exchangeable and in general asymptotically correct. Here, $c_{w,p}(1-\alpha_2)$ denotes the $(1-\alpha_2)$-quantile of the permutation distribution of $T_w$.

{\color{black} Finally, we join the test decisions of the suggested procedures as explained in Section \ref{int} to obtain two different MCTs for testing $H_0^\mu$ against $H_1^\mu:\{\mu_1>\mu_2\}$:
$\varphi_{tw} = \varphi_{t}\cdot \varphi_{w}$ and $\varphi^\mu = \varphi_{t,p}\cdot \varphi_{w,p}$. Both are asymptotic level $\alpha$ tests under $H_0^\mu$ while 
$\varphi^\mu$ may even be exact if the defining permutation tests are.}
The two sided test can be obtained similarly.\\


\section{Nonparametric Test Procedures} \label{nonparametric}
{\color{black}
More general as in Section~\ref{semi} we now consider a nonparametric set-up, i.e. we assume Model (\ref{model: missing}) with arbitrary unknown marginal distribution functions $F_i$ for component $i=1,2$. The only requirement is that $F_i$ is no one point distribution. Here, the main aim is to construct a MCT for the null hypothesis  $H_0^p: \{p= 1/2 \}$ defined in the Wilcoxon-Mann-Whitney (WMW) effect $p=\int F_1 dF_2$. However, the adaptability of the idea will first be exemplified using two other testing problems.

\subsection{Wilcoxon rank sum procedures}\label{WMW}

In this subsection, we specialize to the classical model of the Wilcoxon signed rank and the WMW test. That is we assume a shift model with $F_1(x)=F_2(x-\delta)$ for some $\delta\geq0$. The corresponding null hypothesis of interest is given by $H_0^{\delta}:\{\delta=0\}$ which may be tested against the one-sided alternative $\{\delta>0\}$. In this set-up the unpaired two-sample problem defined via the 
incompletely observed pairs $\vX^{(i)}$ may be tested by means of the WMW test. 
Denoting the mid-rank of $X^{(i)}_{j \ell }$ ($j=1,2, \ \ell=1,\dots,n_\ell$) under all $M=n_1+n_2$ observations in $\vX^{(i)}$ 
by $R_{j \ell}$ it is given by $\varphi_{WMW} = 1\{W>z_{1-\alpha_2}\}$. Here 
$W = \sqrt{n_1 n_2}(\overline{R}_{2 \cdot} - \overline{R}_{1 \cdot})/\sqrt{\widehat{\sigma}_0M}
$ for the rank means $\overline{R}_{j \cdot} = n_j^{-1}\sum_{\ell=1}^{n_j} R_{j \ell}$ and the 
tie-respecting variance estimator 
$\widehat{\sigma}_0^2 = \frac1{M-1} \sum_{j=1}^2 \sum_{\ell=1}^{n_j} \left( R_{j \ell} - ({M+1})/2 \right)^2$. In case of continuous distribution functions $F_1$ and $F_2$ the latter simplifies to the well known factor $\widehat{\sigma}_0^2 = M(M+1)/12$. Moreover, for smaller sample sizes, one usually utilizes the permutation version 
of $\varphi_{WMW}$ since data in $\vX^{(i)}$ is exchangeable under $H_0^\delta$.
For the completely observed pairs  $\vX^{(c)}$ textbooks usually recommend the application of the 
 Wilcoxon signed rank test $\varphi_{WSR}=1\{W^+>z_{1-\alpha_1}\}$. It is based on calculating the ranks $R^+_g$ of the absolute values of 
 $\Delta_g = X_{1g}^{(c)} - X_{2g}^{(c)}$ which are not equal to zero. 
 Denoting the sum of the ranks $R_g^+$ that belong to positive (negative) values of $\Delta_g$ as $R^+$ ($R^-$) the test statistic is given by 
 $W^+ = (R^+ - R^-)/ \sqrt{\sum_g R_g^{+ \, 2}}$. Again, a permutation version is usually performed for smaller sample sizes (assuming $0$-symmetry of $\Delta_g$). 
 The resulting MCT is given by $\varphi_W = \varphi_{WSR}\cdot \varphi_{WSW}$. 
 Since $\varphi_{WSR}$ is based on the differences $\Delta_g$, the resulting test is, however, no test for ordinal data. 
 This unfortunate restriction is inherited from the Wilcoxon signed rank test and relaxed in the next section.

\subsection{Test procedures for $H_0^F:F_1=F_2$}\label{H0F}

We now consider the null hypothesis $H_0^F:\{F_1=F_2\}$ in the more general model described at the onset of Section~\ref{nonparametric}. For this testing problem the 
 WMW-test as defined in Subsection~\ref{WMW} remains valid for comparing the incompletely observed pairs, see e.g. \cite{hollander2013nonparametric}. However, in order to obtain a procedure that is valid for ordinal data, we consider the paired rank test proposed by \cite{munzel1999nonparametric}. To describe the procedure, let $R_{j g}^{(c)}$ be the rank of $X_{j g}^{(c)}$ under all $2n_c$ pooled observations in 
$\vX^{(c)}$. Set $D_g = R_{2 g}^{(c)} - R_{1 g}^{(c)}, g=1,\dots,n_c$ and denote their mean and empirical variance as $\overline{D}_\cdot$ and $\widehat{\sigma}_D^2$, respectively. Then the Munzel-statistic is given by $T_{M,F} = \sqrt{n_c} \ \overline{D}_\cdot/ {\widehat{\sigma}_D}$.
}
 Under the null hypothesis $H_0^F$, this statistic $T_{M,F}$ follows asymptotically (as $n\to\infty$) the standard normal distribution. Therefore, an asymptotic level $\alpha_1-$test for the one sided case is given by $\varphi_{M,F}=1\{T_{M,F} \geq z_{1-\alpha_1}\}$.
 {\color{black}For smaller sample  sizes \cite{munzel1999nonparametric}
suggested to substitute the standard normal quantile $z_{1-\alpha_1}$ with a $t_{n_c-1, 1-\alpha_1}$ quantile. Finally, an asymptotic level $\alpha$ MCT for $H_0^F$ is given by $\varphi_{F}= \varphi_{M,F}\cdot \varphi_{WMW}$. Different to $\varphi_W$, the MCT $\varphi_F$ is also applicable for ordinal data.\\
}

{\color{black} Since {\it 'effect sizes are the most important outcome of empirical studies'} \citep{lakens2013calculating} 
we close this section with a detailed study of the WMW effect.
}

\subsection{Test procedures for $H_0^p:p=1/2$}
{\color{black} 
To entail a parameter for describing differences between distributions we consider the WMW-effect
\bqan \label{treatment effect}
p=\int F_1dF_2=P(X_{11}^{(c)}<X_{22}^{(c)})+1/2P(X_{11}^{(c)}=X_{22}^{(c)}),
\eqan 
also known as treatment effect or relative marginal treatment effect \citep{fligner1981robust,  brunner199619, munzel1999nonparametric, brunner2000nonparametric}. 
Its interpretation is rather simple: 
If $p>1/2$, the observations from $F_2$ tend to be larger than those from $F_1$, and vice versa if $p<1/2$. Furthermore, $p=1/2$ corresponds to 
the case of no treatment effect. Our aim is to construct a test for the so-called nonparametric Behrens-Fisher problem 
$H_0^p:\{p=1/2\}$ in the general nonparametric Model (\ref{model: missing}). To this end, we mainly have to change the variance estimators of the test statistics presented in Section~\ref{H0F}. In particular, we infer the complete cases $\vX^{(c)}$ by means of the 
\cite{munzel1999nonparametric} statistic for the paired Behrens-Fisher problem given by

 \bqan \label{Munzel test}
T_{M}= T_{M}(\vX^{(c)}) = \sqrt{n_c} {\overline{D}_{\cdot}}/{s_c},
\eqan
where 
$s_c^2=(n_c-1)^{-1}\sum_{g=1}^{n_c} (Z_g-\overline{Z}_\cdot)^2$ is the empirical variance of 
$Z_g=(R_{2g}^{(c)}-R_{2g}^{(c,2)}-(R_{1g}^{(c)}-R_{1gk}^{(c,1)}))/n_c$. Here, 
$R_{jg}^{(c,j)}$ is the internal rank of $X_{jg}^{(c)}$ among the $n_c$ observations $X_{j1}^{(c)},....,X_{jn_c}^{(c)}$ from component $j = 1,2.$. 
Under the null hypothesis $H_0^p$, the test statistic $T_{M}$ follows asymptotically (as $n\to\infty$) the standard normal distribution. 
Therefore, \cite{munzel1999nonparametric} proposed 
$\varphi_{M,1}=1\{T_M \geq z_{1-\alpha_1}\}$ and $\varphi_{M}=1\{T_M\leq -z_{1-\alpha_1/2}\}+1\{T_M\geq z_{1-\alpha_1/2}\}$ as asymptotic level $\alpha_1$ procedures for testing one-sided ($H_{1,>}^p:\{p>1/2\}$) and two-sided alternatives ($H_1^p:\{p\neq1/2\}$), respectively. 
 However, for small to moderate sample sizes, previous simulation studies {\color{black} of \cite{konietschke2012studentized}} showed that these tests do not control the type-I error rate constantly under a range of correlation values. Therefore, to enhance it's finite sample performance, 
they recommended the use of a studentized permutation version which is based on randomly changing the components for each pair. It is given as 
{\color{black}$\varphi_{KP}^{(c)} =1\{T_M\leq z_{\alpha_1/2}^{\tau}\}+1\{T_M\geq z_{1-\alpha_1/2}^{\tau}\}$, where $z_{1-\alpha_1/2}^{\tau}$ denotes the $(1-\alpha_1/2)$-quantile} of $T_M$'s permutation distribution. As shown in \cite{konietschke2012studentized} this test is asymptotically exact and may even be finitely exact in case of invariance. \\
}
{\color{black} Considering the incompletely observed vectors $\vX^{(i)}$ we encounter the situation of the 
classical nonparametric Behrens-Fisher problem considered in \cite{brunner2000nonparametric}. They proposed to use the test statistic
\bqan \label{NB test}
T_{BM}= T_{BM}(\vX^{(i)}) = \sqrt{\frac{n_1n_2}{M}} \ \frac{\overline{R}_{2\cdot}-\overline{R}_{1\cdot}}{s_i}.
\eqan
Therein, $s_i^2=M(\frac{1}{n_2}\hat{\sigma}_1^2+\frac{1}{n_1}\hat{\sigma}_2^2)$ for 
 $\hat{\sigma}_j^2=\frac{1}{n_j-1}\sum_{k=1}^{n_j}(R_{jk}-R_{jk}^{(i)}-\overline{R}_{j\cdot.}+\frac{n_j+1}{2})^2$ and 
 $R_{jk}^{(i)}$ is the internal rank of $X_{jk}^{(i)}$ among the $n_j$ observations $X_{j1}^{(i)},....,X_{jn_j}^{(i)}$ for $j=1,2$. Again $T_{BM}$ is asymptotically standard normal under $H_0^p$ and \cite{brunner2000nonparametric} even proposed a $t$-approximation for smaller sample sizes. Moreover, the latter was slighlty enhanced by a studentized permutation approach leading to {\color{black}the two-sided} permuted Brunner-Munzel test 
{\color{black} $\varphi_{NB} =1\{T_{BM}\leq z_{\alpha_2/2}^{\pi}\}+1\{T_{BM}\geq z_{1-\alpha_2/2}^{\pi}\}$, , where $z_{1-\alpha_2/2}^{\pi}$ is the corresponding $(1-\alpha_2/2)$- permutation quantile}, see \cite{neubert2007studentized, pauly2016permutation}.
 The MCT for testing the hypothesis $H_0^p$ is thus given by $\varphi^{(p)}= \varphi_{KP}\cdot \varphi_{NB}$.\\
 
 We finally like to emphasize that MCTs for hypotheses formulated in terms of effects can also be inverted into confidence intervals. Those are obtained by unifying the confidence intervals for the underlying paired and unpaired two-sample problems. 


\section{Comparison with existing procedures} \label{alt}

In this section we review existing methods for paired designs with missing values. 
Due to abundance of different methods we only consider the semiparametric ($H_0^\mu$) and one of the nonparametric hypotheses. Due to its broad applicability and easy interpretation of the underlying effect $H_0^p$ got the vote over 
$H_0^\delta$ and $H_0^F$. 
Major assessments were availability of meaningful simulation results from existing literature together with the procedures' validity under the rather general model assumptions. Most tests which do not fulfill both criteria were excluded from the subsequent discussion.
}

\subsection{Test procedures for $H_0^\mu:\mu_1=\mu_2$}
For the semiparametric case described in Section~\ref{semi}, few procedures for mean-based inference in partially observed pairs exist. In particular, \cite{lin1974difference} and \cite{ekbohm1976comparing} developed tests based on simple mean difference estimators which can be transferred to the present situation, see also \cite{uddin2017testing}. Extensive simulation results by \citet{bhoj1991testing, kim2005statistical, xu2012accurate, samawi2014notes} and \cite{amropauly2017}, however, recommend weighted combinations of test statistics for the corresponding paired and unpaired two-sample problems. 
Thereof, we choose the permutation procedure in the weighted test statistic
\bqan \label{T}
T_{ML} = \sqrt{a} T_{t}(\vX^{(c)}) + \sqrt{1-a} T_{w}(\vX^{(i)})
\eqan
proposed by \cite{amropauly2017} which was found to be robust against heteroscedasticity and skewed distributions. 
Here, the weight was chosen as {$a=2nc/(n+nc)$} as recommended in that paper and permutation was achieved by flipping components in the completely observed variables $\vX^{(c)}$ and randomly shuffling the group status in the incompletely observed cases $\vX^{(i)}$. 
We will compare this procedure with the MCTs $\varphi_{tw}$ and $\varphi_{\mu}$ based on paired- and Welch-t-tests with $t$- ($\varphi_{tw}$) or permutation quantiles ($\varphi_{\mu}$) as critical values in the previous section.

\subsection{Test procedures for $H_0^p:p=1/2$}

Different to the nonparametric hypothesis $H_0^F$ (see e.g. \cite{fong2017rank}) there is no abundance of methods 
testing the null hypothesis $H_0^p:\{p=1/2\}$ in matched pairs with missing data. In particular, a thorough literature review 
only left the nonparametric ranking procedure by \cite{konietschke2012ranking} as adequate competitor since it has shown the best performance in comparisons to other rank-based procedures \citep{konietschke2012ranking}. It is based on the test statistic
\bqan \label{KHLB}
T_{KHLB}^{p}= \sqrt{n}\ \frac{\hat{p}_\theta-1/2}{\tilde{s}_\theta},
\eqan 
where $\hat{p}_\theta=\frac{1}{N}\bigg(\theta_2\overline{R}_{2.}^{(c)}-\theta_1\overline{R}_{1.}^{(c)}
+(1-\theta_2)\overline{R}_{2.}^{(i)}-(1-\theta_1)\overline{R}_{1.}^{(i)}
\bigg)+\frac{1}{2}$ is a weighted estimator of the WMW-effect $p$ for the weights $\theta_j=n_c/(n_c+n_j)$, $j=1,2$. The explicit formula for the rank-based standard deviation $\tilde{s}_\theta$ can be found in 
{\cite{konietschke2012ranking}.} 
Under $H_0^p$, the statistic $T_{KHLB}^{p}$ has asymptotically a standard normal distribution as $n_c+n_j\rightarrow \infty, j=1,2$. Based on their simulations, an additional $t$-approximation was recommended to ensure better control of the type-I error rate in case of small to moderate sample sizes. The corresponding procedures are denoted as $T_{KHLB}(normal)$ and $T_{KHLB}(t-App)$, respectively.


\begin{figure}
	\begin{center}
	\includegraphics[width=1\textwidth]{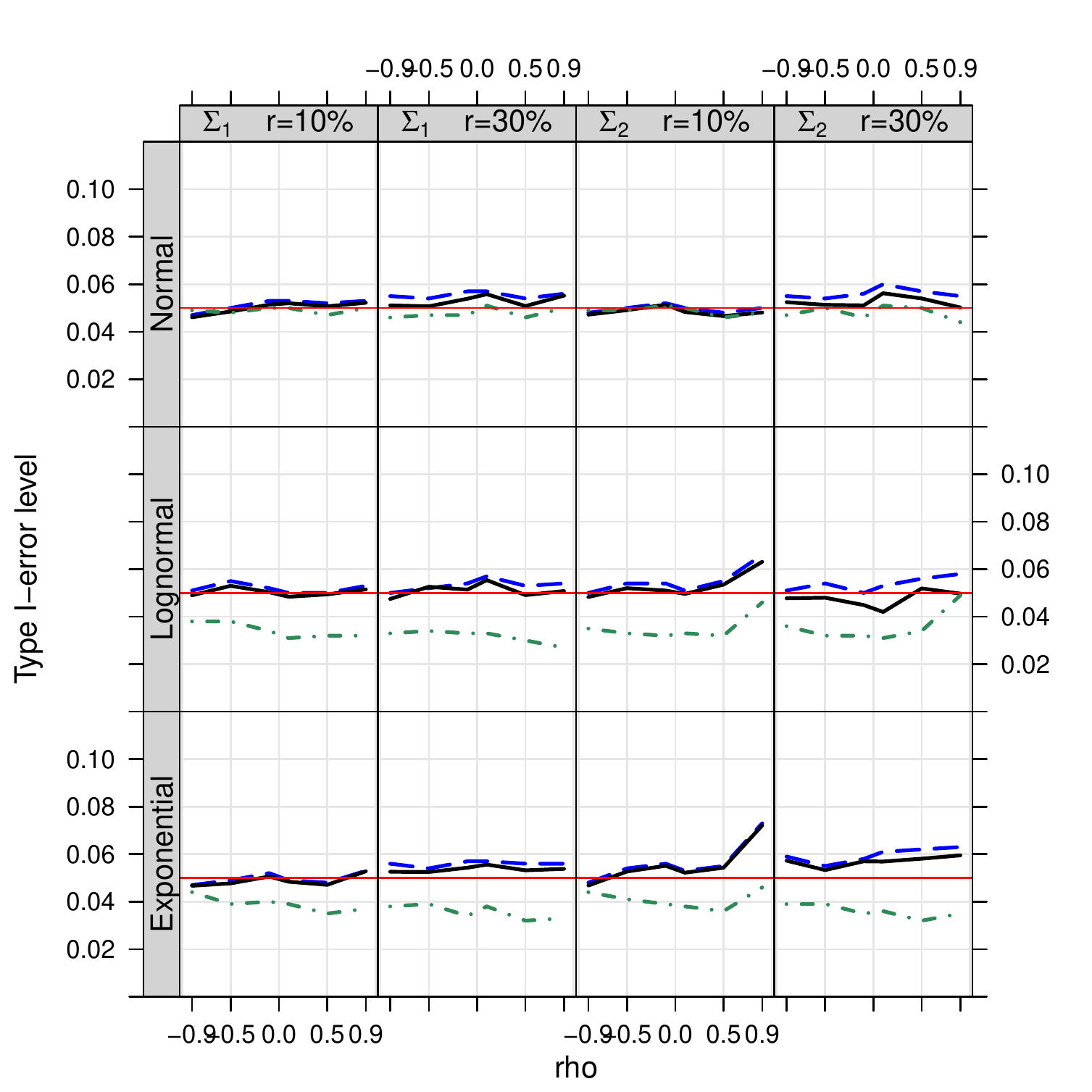}
		\caption{Simulation results for the type-I error level for $\varphi^\mu$ $(\textendash\textendash\textendash)$, $T_{ML}$ ({\color{black}$--$}), and 
			$T_{tw}$ ({\color{ao(english)}$-\cdot-$}) for different distributions under varying correlation factors $\vSigma=\vSigma_1$ (left) and $\vSigma=\vSigma_1$ (right)  for testing $H_0^\mu:\{\mu_1=\mu_2\}$ with $n=10$ (left) and $n=20$ (right) and different missing percentages $r\in\{10\%, 30\%\}$ under the MCAR framework.}
		\label{fig:multiH0m}
	\end{center}	
\end{figure}


\section{Monte-Carlo simulations} \label{sim}
\begin{figure}
	\begin{center}

\includegraphics[width=1\textwidth]{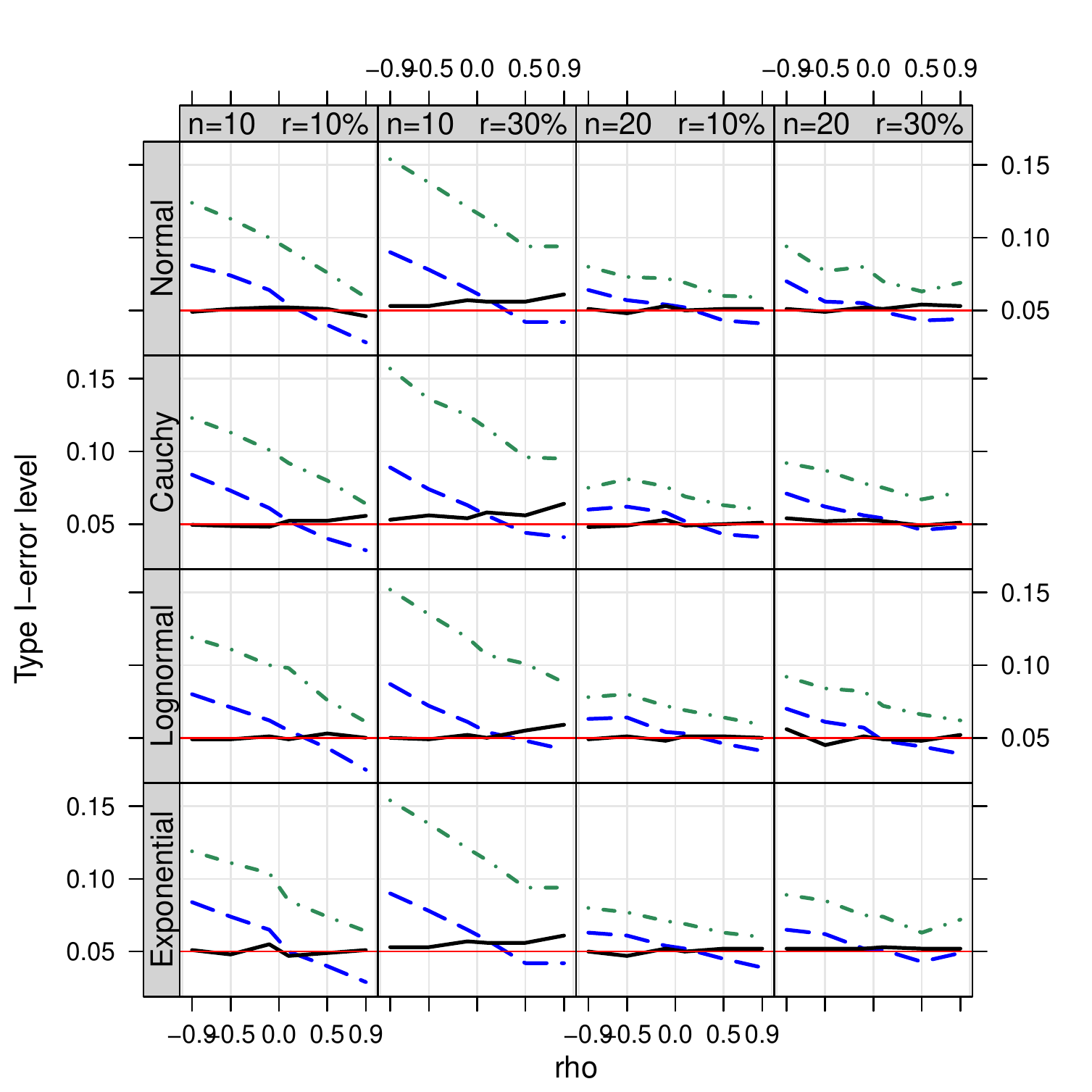}
	
		\caption{Simulation results for the type-I error level for $\varphi^{(p)}$ $(\textendash\textendash\textendash)$, $T_{KHLB}(t-App)$ ({\color{black}$--$}), 
			and $T_{KHLB}(normal)$ ({\color{ao(english)}$-\cdot-$}) for different distributions under varying correlation factors {\color{black}($\vSigma=\vSigma_1$)} for testing $H_0^p:\{p=1/2\}$ with $n=10$ (left) and $n=20$ (right) and different missing percentages $r\in\{10\%, 30\%\}$ under the MCAR framework.}
		\label{fig:multi}
	\end{center}	
\end{figure}

\begin{figure}
	\begin{center}
		
		\includegraphics[width=1\textwidth]{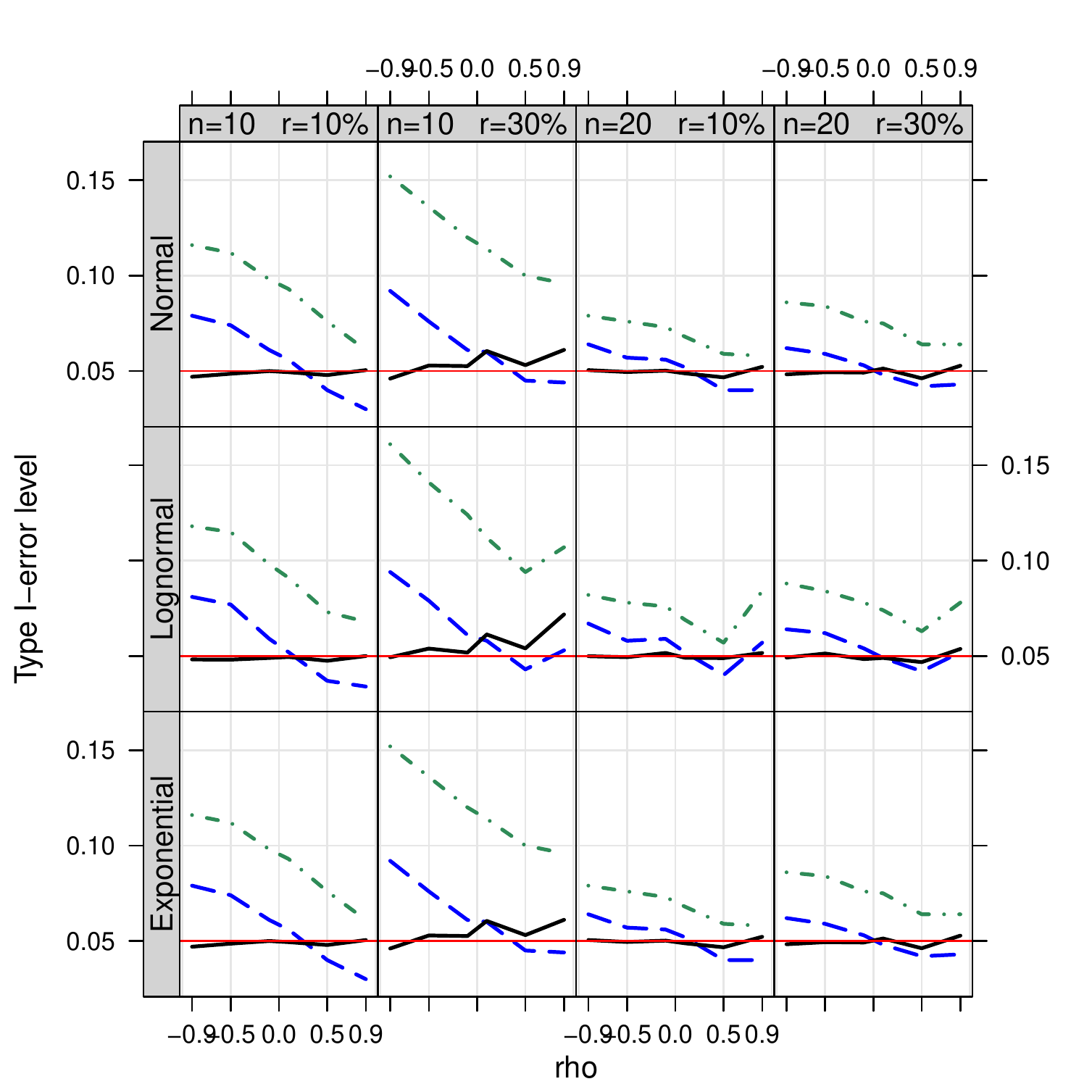}
		
		\caption{Simulation results for the type-I error level for $\varphi^{(p)}$ $(\textendash\textendash\textendash)$, $T_{KHLB}(t-App)$ ({\color{blue}$--$}), 
			and $T_{KHLB}(normal)$ ({\color{ao(english)}$-\cdot-$}) for different distributions under varying correlation factors {\color{black}($\vSigma=\vSigma_2$)} for testing $H_0^p:\{p=1/2\}$ with $n=10$ (left) and $n=20$ (right) and different missing percentages $r\in\{10\%, 30\%\}$ under the MCAR framework.}
		\label{fig:multivar2}
	\end{center}	
\end{figure}

A simulation study was conducted to investigate the finite sample behavior of our suggested MCTs and the selected procedures from Section \ref{alt} for testing $H_0^\mu$ and $H_0^p$, respectively. {\color{black} All tests are performed as two sided-tests}. 
Major assessment criterion was their type-I error rate control at level $\alpha=0.05$. 
 The complete pairs were generated by
 \begin{center} 
	$\begin{pmatrix}
	X_{1j} \\ X_{2j} 
	\end{pmatrix}=  \begin{pmatrix}
	\mu_1\\ \mu_2
	\end{pmatrix} + \vSigma^{1/2} \begin{pmatrix}\varepsilon_{1j} \\\varepsilon_{2j} \end{pmatrix}, \quad j=1,\dots,n$
\end{center}
 with standardized i.i.d. errors $\varepsilon_{ij}, i=1,2, j=1,\dots,n$. In particular, we simulated 
 
 symmetric as well as skewed $\varepsilon_{ij}$: Standardized normal, exponential, log normal, and plain Cauchy distribution. For $\vSigma$ we considerd the choices  
\begin{center}
	$\vSigma_1 = \begin{pmatrix} 1&\rho\\ \rho&1 \end{pmatrix}$   \hspace{2pt}\text{and}\hspace{2pt}
{\color{black}	$\vSigma_2 = \begin{pmatrix} 1&\sqrt{2}\rho\\ \sqrt{2}\rho& 2\end{pmatrix}$ }
\end{center}
with varying correlation factor $\rho \in (-1,1)$, representing a homoscedastic and a heteroscedastic setting, respectively. We note, that $\rho$ corresponds to, but is in general not, the correlation \citep{leonov2017correlated}. The sample sizes were chosen as $n \in \{10, 20\}$ and missingness was included within a MCAR framework, i.e.
observations from Model \eqref{model: missing} were simulated as
$\vX_k=(\delta_{1k}X_{1k},\delta_{2k}X_{2k}), k=1,...,n$ for  Bernoulli distributed $\delta_{ik}\sim B(r)$. 
The missing probability was chosen to $r=10\%, 30\%$ and the independent tests for the corresponding paired and unpaired testing problem 
were calculated at significance level $\alpha_1=\alpha_2=\sqrt{\alpha}$. We note that other choices of 
$\alpha_1$ and $\alpha_2$ lead to similar results in our scenarios (results not shown). 
All computations were done using the R computing environment \citep{Rlanguage}, version 3.2.2 (on $10,000$ simulation runs and $B=1,000$ permutations runs).\\

The type-I error level $(\alpha=0.05)$ simulation results for testing $H_0^\mu$ and $H_0^p$ are summarized in {\color{black}Figures \ref{fig:multiH0m} - \ref{fig:multivar2} and described below.}

{\it Semiparametric mean-based case:} 
It can be readily seen from Figure~\ref{fig:multiH0m} that the MCT based on the classical paired-$t$- and Welch-$t$-tests tends to rather conservative conclusions when data is not normally distributed. In contrast the two other methods based on permutation quantiles both 
control the nominal type-$I$ error rate accurately. In particular, they were only slightly affected by the amount of missingness. Only in case of skewed distributions appearing together with large positive correlation factors ($\rho=0.9$) 
a slightly liberal behviour can be observed.  

{\it Nonparametric case:} The simulation results for testing $H_0^p$ show that the new MCT $\varphi^{(p)}$ greatly improves upon the existing methods by \cite{konietschke2012ranking}. It can be seen from  {\color{black} Figure~\ref{fig:multi} (homoscedastic settings) and Figure~\ref{fig:multivar2} (heteroscedastic settings) } that the two tests based on the statistic $T^p_{KHLB}$ in \eqref{KHLB} tend to be very sensitive to the dependency structure in the data. In particular, {\color{black}under homoscedasticity as well as heteroscedasticity} both procedures exhibit a liberal behaviour for negative correlation factors $\rho$ and a conservative one for positive $\rho$ which is more pronounced for $T_{KHLB}^p(normal)$. Contrary, the MCT $\varphi^{(p)}$ exhibits a fairly good type-$I$ error rate control for all considered scenarios. Only in case of large correlation factors ($\rho\approx 0.9$) and higher missing rate ($r=30\%$) a slightly liberality can be observed for small sample sizes ($n=10$).

\section{Analysis of the data example}\label{example}
{\color{black}In this Section, we reconsider the clinical migraine study by \cite{kostecki1999pine}. The study has been performed to investigate four consecutive sessions of a non-drug headache program. 
Hereby, the headache severity level of $135$ migraine patients was
 measured repeatedly on a daily basis by a scale ranging from $0$ to $20$. The lower the score, the better the clinical record. The goal of this study was to monitor the changes of headache severity level of the patients during the four sessions and to determine the efficiency. However, many patients missed the record measurements for some sessions 
 leading to a large amount of missing observations. \cite{gao2007nonparametric} investigated the missing data mechanism by performing correlation tests for the missing proportions versus the average headache severity scores across the sessions. She concluded that the MCAR assumption is reasonable for analyzing the data.\\
 For reformatting the data into a matched pairs design we only use the clinical records of the patients in their first and third session. A close look to our data shows that from the total of $132$ patients who attained Session 1 and Session 3, only $n_c=82$ patients were measured twice, $n_1=44$ patients were only observed on Session 1 $n_2=6$ patients were only scored in Session 3. The data were analyzed by all our considered testing methods $T_P, T_{KHLB}(t-app)$ and $T_{KHLB}(normal)$ for the null hypothesis $H_0^p:\{p=1/2\}$. The results are summarized in Table~\ref{pvalue}. It can be seen that all tests indicate a significant difference between
 the headache severity levels of the two sessions (Two sided p=value $<0.01$). Moreover, based on the results of the right sided test, we conclude that the clinical outcome of the migraine patients significantly improves after two sessions (One sided p=value $<0.01$).}

\begin{table}
	\tbl{P-values of  the migraine patients Data} 
	{\begin{tabular}{l c c } 
			\hline\hline 
			Method & One-sided p-value & Two-sided p-value \\ [0.5ex] 
			\hline \hline 
			MCT $\varphi^{(p)}$ & 0.0033 & 0.0066\\ 
			$T_{KHLB}(t-app)$ & 0.0027 & 0.0053\\
			$T_{KHLB}(normal)$ & 0.0023 & 0.0045\\
			\hline \hline 
		\end{tabular}}
		\label{pvalue} 
	\end{table}

\section{Discussion and Outlook} \label{dis}

In this paper, we investigated statistical inference methods for partially observed pairs. Under the MCAR framework, a  {\it multiplication-combination test (MCT)} procedure that is flexibly applicable to test hypotheses formulated in terms of various effect sizes (statistical functionals) has been proposed. It is based on combining the results of two independent tests for the related paired and unpaired two-sample problem, say $\varphi_1$ and $\varphi_2$. 
The key idea is to carry out each test $\varphi_i$ at larger significance levels $\alpha_i$ ($i=1,2$) such that $\alpha_1\alpha_2=\alpha$ (e.g. $=0.05$). Since larger percentiles ($\alpha_i$) can usually be estimated more accurately than smaller ones ($\alpha$) we expected a better type-$I$ error control in case of small to moderate sample sizes. The idea has been applied in different statistical models including semi- as well as purely nonparametric models. In particular, for testing the Behrens-Fisher problem $H_0^p:\{p=1/2\}$, a MCT based upon permutation versions of the \cite{munzel1999nonparametric} (paired) and \cite{brunner2000nonparametric} (unpaired) considerably improved existing methods. This is especially astonishing since \cite{fong2017rank} have recently concluded that 
{\it 'a potential obstacle in doing so is the inflated type 1 error rates ($\ldots$) which occur because the complexity of the variance formulae requires bigger sample size to justify large sample approximation'}. 

Future research will be concerned with less stringent missing mechanisms (as MAR), different effect sizes (e.g., the $\log$-odds-ratio) and more complex designs (e.g., more groups or time points).

\section*{Acknowledgement}
This work was supported by the German Academic Exchange Service (DAAD) under the project: Research Grants - Doctoral Programmes in Germany, 2015/16 (No. 57129429). 
\bibliographystyle{plainnat}
\bibliography{References2}

\end{document}